\documentclass[11pt]{article}
\usepackage{a4wide}
\usepackage{amsmath}

\usepackage{amsfonts}
\usepackage{amssymb}
\usepackage{euscript}
\newenvironment{proof}{\noindent {\it Proof.~~}\ }{\  \rule{1mm}{2mm}\medskip}

\newenvironment{proofof}[2]{\noindent {\it Proof of #1}~#2: \
}{~\rule{1mm}{2mm}\medskip}
\newtheorem{theorem}{Theorem}
\newtheorem{lemma}[theorem]{Lemma}
\newtheorem{corollary}[theorem]{Corollary}
\newtheorem{proposition}[theorem]{Proposition}

\newcommand{\subgp}[1]{\langle{#1}\rangle}

\begin{document}
\title{
  Hyper-atoms and the Kemperman's critical pair Theory}

\author{ Yahya O. Hamidoune\thanks{Universit\'e Pierre et Marie Curie,    Paris {\tt hamidoune@math.jussieu.fr} }
}
%soumis Combinatorica, Receipt, August,27-2007
\date{August, 15-2007}
\maketitle

\begin{abstract}
In the present work, we  introduce   the notion of a hyper-atom and
prove their main structure theorem. We then apply the global
isoperimetric methodology to give a new proof for
 Kemperman's
structure Theory and a  slight improvement.

\end{abstract}

\section{Introduction}

A basic tool in Additive Number Theory is the following
generalization of the Cauchy-Davenport Theorem
\cite{cauchy,davenport} due to Kneser:

\begin{theorem}[Kneser \cite{natlivre,tv}]\label{kneser}
Let $G$ be an abelian group and let $A, B\subset G$ be finite
subsets of $G$.  Then $|A+B|\ge |A+H|+|B+H|-|H|$, where $H$ is the
period of $A+B$.
\end{theorem}

The description for the subsets with $|A+B|= |A|+|B|-1$  needs some
terminology:

A decomposition
 $A=A_0\cup A_1$ is said to be a {\em $H$--quasi-periodic
 decomposition} if
 $A_0+H=A_0$ and $A_1$ is contained in some $H$--coset.
 Let $A,B\subset G$.

A pair $\{A,B\}$ will be called an  {\em elementary pair} if one of
the following conditions holds:
\begin{itemize}

\item [(SP1)] there is $d\in G$, with order $\ge |A|+|B|-1$ such that $A$ and $B$ are arithmetic progressions with difference
$d$,
\item [(SP2)] $\min(|A|,|B|)=1$,
\item [(SP3)] $A$ is aperiodic and there is a finite subgroup  $H$ and $g\in G$ such that $A,B$ are contained in some
$H$--cosets and  $g-B=H\setminus A,$ and for all $c$,
$|(c-A)\cap B|\neq 1$,
\item [(SP4)] there is a subgroup  $H$ with a prime order such that $A,B$ are contained in some $H$-cosets and $|A|+|B|=|H|+1
$, and  moreover there is a unique $c\in G$ such that
$|(c-A)\cap B|=1.$

 %Moreover  $|(c-A)\cap B|\neq 1,$ for all $c$.
\end{itemize}

An elementary pair satisfying one of the conditions (SP1), (SP2) or
(SP3) will be called a {\em strict elementary pair}.

Notice that the condition "with a prime order" in SP4 is not present
  in Kemperman's formulation. Hence the class of elementary pairs in the sense of Kemperman is
larger than our class. This will produce a slightly more precise
result than the result proved by Kemerman:

\begin{theorem}\label{Kemperman} \cite{kempacta}
Let $A,B$ be finite subsets of an abelian group $G$ with $|G|\ge 2.$

Then  the following conditions are equivalent:

\begin{itemize}
\item [(I)] $|A+B|= |A|+|B|-1,$ and moreover $|(c-A)\cap B|=1$ for some
$c$ if $A+B$ is periodic.

\item [(II)]  There is a nonzero subgroup $H$   and  $H$-quasi-periodic decompositions
 $A=A_0\cup A_1$ and $B=B_0\cup B_1$ such that $(A_1,B_1)$ is an
 elementary pair and $|(\phi(a_1+b_1)-\phi(A))\cap \phi (B)|=1$,
 where $\phi : G\mapsto G/H$  is the canonical morphism and
 $a_1\in A_1$ and $b_1\in B_1$.
 \end{itemize}

\end{theorem}

 The redundant condition $|\phi(A+B)|=|\phi(A)|+|\phi(B)|-1$ present in
Kemperman's formulation was omitted  since it is a consequence of
the condition "$|(\phi(a_1+b_1)-\phi(A))\cap \phi (B)|=1$" by
Scherck's Theorem \ref{scherk}. The original and unique previously
known proof of Kemperman's result uses the additive local
transformations introduced by Cauchy and Davenport
\cite{cauchy,davenport}.

Recently the author introduced the isoperimetric method allowing to
derive additive inequalities from global properties of the fragments
and atoms (subsets where the objective function $|A+B|-|A|$ achieves
its  non trivial minimal value).

This method can be applied to abstract graphs and non abelian groups
and have implications that could not be derived  using the local
transformations. However in the abelian case, it was not clear how
to derive the Kneser-Kemperman Theory using this  method.

Very recently Balandraud introduced some isoperimetric objects and
proposed a proof, requiring several pages, of Kneser's Theorem using
as a first step our result that the $1$-atom containing $0$ is a
subgroup.

On the other side, alternative proofs for results proved first using
the isoperimetric method, based on  Kemperman's Theory as a main
tool, were obtained by Grynkiewicz in \cite{davdecomp} and Lev
\cite{levkemp}.

In the present work, we  introduce   the notion of a hyper-atom and
prove the main structure theorem for hyper-atoms. We then apply the
global isoperimetric methodology introduced in \cite{hglobal}
 to give a new proof for
 Kemperman's
structure Theory with a slight improvement.  The methods introduced
in the present work allow quite likely much more complicated
descriptions for subsets $A,B$ with $|A+B|=|A|+|B|+m$, with small
some small values of $m\ge 0$. We made the calculations for $m=0$,
obtaining a new proof of a recent result due Grynkiewicz in
\cite{davkem}, that extends to all abelian groups  a result proved
by R{\o}dseth and the author \cite{hrodseth2}. However we shall
limit ourselves to  Kemperman's Theory in order to illustrate the
method in a relatively simple context.

\section{Terminology and preliminaries}

Let  $ A,B$ be subsets of $ G $. The subgroup generated by  $A$ will
be denoted by $\subgp{A}$. The {\em Minkowski sum} is defined as
$$A+B=\{x+y \ : \ x\in A\  \mbox{and}\ y\in
  B\}.$$

Let  $H$ be a subgroup. A partition  $A=\bigcup \limits_{i\in I}
A_i,$ where $A_i$ is the nonempty intersection of some $H$--coset
with $A$ will be called
 a $H$--{\em decomposition}
of $A$.

%Let $X$ be a subset of a group $G$. We write $\Pi(X)=\{x : X+x=X\}.$

For an element $x\in G$, we write $r_{A,B}(x)=|(x-B)\cap A|$. Notice
that $r_{A,B}(x)$ is the number of distinct representations of $x$
as a sum of an element of $A$ and an element of $B$.

We use the following well known and easy fact:
\begin{lemma}\cite{natlivre}
Let $G$ be a finite group and let
 $A,B$  be  subsets
 such that $|A|+|B|\ge |G|+t$.
 Then  $r_{A,B}(x)\ge t$.

\label{prehistorical}
 \end{lemma}

We shall use the following result:

\begin{theorem}\label{scherk}(Scherk)\cite{scherck}.
  Let $X$ and $Y$ be nonempty finite subsets of an abelian group $G$. If there is
an element   $c$   of  $G$ such
  that $|X\cap(c-Y)| = 1$, then
    $|X + Y| \geq |X| + |Y| - 1.$
\end{theorem}

Scherck's Theorem follows easily from Kneser's Theorem, c.f.
\cite{gerlodinger}. We give in the appendix a short direct proof for
this result.

We need Vosper's Theorem:

\begin{theorem}\label{vosperth}
Let  $A,B$ be subsets of a group $G$ with a prime order such that
 $|A|,|B|\ge 2$ and $|A+B|= |A|+|B|-1\le |G|-2.$ Then  $A,B$  are  arithmetic progressions with the same difference.

\end{theorem}
As showed in \cite{halgebra,hiso2007}, this result follows in few
lines from the intersection property of the $2$--atoms.

%In the case $|A|+|B|= |G|+1$, we still having some structure under
%an additional condition:

%\begin{lemma}\label{vosperiv}
%Let  $A,B$ be subsets of a group $G$ with a prime order such that $
%|A|+|B|+1= |G|.$ Assume  moreover there is a unique $c\in G$ such
%that $|\{ c : |(c-A)\cap B|=1\}|\ge 2.$

%Then  either $\min (|A|,|B|)=1$ or $A,B$  are  arithmetic
%progressions with the same difference.

%\end{lemma}

%This lemma will be proved in the appendix.

 Let $V$ be a set and let $E \subset V\times V$.  The relation
$\Gamma = (V,E)$ will be called  a  {\em graph}.
 The elements of  $V$
will be called {\em points}. The graph $\Gamma$ is said to be {\em
reflexive} if $(x,x)\in E,$ for all $x$. We shall write $$\partial
(X)=\Gamma (X)\setminus X.$$

%An  {\em elementary path} of $\Gamma$ from $x_1$ to $x_k$ is a
%sequence $\mu =[x_1, \cdots ,x_k]$ of pairwise distinct points
%(where $k\ge 1$) such that $(x_i,x_{i+1})\in E$, for all $1\le i \le
%k-1.$ The set of points of $\mu$ is by definition $P(\mu)=\{x_1,
%\cdots, x_k\}$.

%A family $\mu _1, \cdots , \mu _k$ of paths from $x$ to $y$ will be
%called openly disjoint if  $P(\mu _i)\cap P( \mu _j)=\{x,y\}$ for
%all $i,j$ with $i\neq j$.

Let $\Gamma =(V,E)$ be a locally finite graph with $|V|\ge 2k-1.$
 The {\em $kth$--connectivity}
of $\Gamma$
 is

$$
\kappa _k (\Gamma )=\min  \{|\partial (X)|\   :  \ \
\infty >|X|\geq k \ {\rm and}\ |X \cup  \Gamma(X)|\le |V|-k\},
$$
where $\min \emptyset =|V|-2k+1$.

Let $G$ be a group, written additively, and let $S$  be a subset of
$G$. The graph $(G,E),$ where $ E=\{ (x,y) : -x+y \ \in S \}$ is
called a {\it Cayley graph}. It will be denoted by $\mbox{Cay}
(G,S)$.

Let $\Gamma =\mbox{Cay} (G,S)$   and  let   $F \subset G $. Clearly
 $\Gamma (F)=F+S $.

A general formalism, including the most recent isoperimetric
terminology  may be found in a the recent paper \cite{hiso2007}.

We recall that  Menger's Theorem which is a basic min-max relation
from Graph Theory  \cite{hglobal,natlivre,tv} has several
implications in number Theory.  We need the following consequence of
Menger's Theorem:
\begin{proposition} \cite{hglobal}{ Let $\Gamma $ be a locally finite   reflexive
graph and let $k$ be a nonnegative integer with
 $k\le \kappa _1$. Let   $X$  a
finite subset of $V$ such that $\min (|V|-|X|, |X|)\ge k.$ There are
pairwise distinct elements
 $x_1, x_2, \cdots, x_{k} \in X$ and pairwise distinct elements
 $y_1, y_2, \cdots, y_{k} \notin X$ such that
 \begin{itemize}
   \item $(x_1,y_1), \cdots , (x_{k}, y_{k})\in E,$
    \item $|X\cup \{y_1, \cdots ,  y_{k}\}|=|X|+k$,
 \end{itemize}

\label{strongip}}
\end{proposition}

We call the property given in Proposition  \ref{strongip}  the {\em
strong isoperimetric property}.

\section{Isoperimetric tools}

The isoperimetric method is usually developed in the context of
graphs. We need in the present work only the special case of Cayley
graphs on abelian groups that we shall identify with group's
subsets.

Throughout all this section, $S$ denotes a finite generating subset
of an abelian group $G$, with $0\in S$.

 For a subset $X$, we put $\partial _S(X)=(X+S)\setminus X$ and $X^S=G\setminus (X+S)$.
%We used the notation $X^\curlywedge$ in \cite{hiso2007} for the
%corresponding notion in arbitrary graphs.
We need the following lemma:

\begin{lemma}\cite{balart,hiso2007}{Let  $X$ be a
 subset of $G$. Then $(X^S)^{-S}+S=X+S$. \label{dualitys}}
\end{lemma}

The last lemma is proved in  Balandraud \cite{balart} and
generalized in \cite{hiso2007}.

 We shall say that a subset
$X$ induces a {\em $k$--separation} if $ |X|\geq k$ and $|X^S|\geq
k$. We shall say that $S$ is $k$--separable if some $X$ induces a
$k$--separation.

Suppose that $|G|\ge 2k-1.$
 The {\em $kth$--connectivity}
of $S$
 is defined  as $\kappa _k (S )=\kappa _k(\mbox{Cay}(G,S))$. By the
 definition we have

$$
\kappa _k (\Gamma )=\min  \{|\partial (X)|\   :  \ \
\infty >|X|\geq k \ {\rm and}\ |X+S|\le |G|-k\},
$$
where $\min \emptyset =|G|-2k+1$.

 A finite subset $X$ of $G$ such that $|X|\ge k$,
$|G\setminus (X+S)|\ge k$ and $|\partial (X)|=\kappa _k(S)$ is
called a {\em $k$--fragment} of $S$. A $k$--fragment with minimum
cardinality is called a {\em $k$--atom}.
%The cardinality of a
%$k$--atom of $S$  will be denoted by $\alpha_k(S)$.

Let $S$ be a non--$k$--separable subset such that $|G|\ge 2k-1$.
Then $G$ is necessarily finite.  In this case, a {\em $k$--fragment}
(resp. {\em $k$--atom})
 is a set with cardinality  $k.$

 These notions, are particular cases  some concepts in
\cite{hcras, Hejcvosp1,halgebra,hactaa,hiso2007}. The reader may
find all basic facts from the isoperimetric method in the recent
paper \cite{hiso2007}.

A $k$--fragment of $-S$ will be  called a {\em negative}
$k$--fragment.
%We use the following notations, where the reference to $S$ could be implicit:
%\begin{itemize}
%\item  $\alpha_{-k}(S ) =\alpha_k(-S),$
 % \item  $\kappa_{-k}(S ) =\kappa_k(-S)$.
%\end{itemize}

Notice that
    $\kappa _k (S)$ is the maximal integer $j$
such that for every finite subset $X\subset G$  with $|X|\geq k$,

\begin{equation}
|X+S|\geq \min \Big(|G|-k+1,|X|+j\Big).
\label{eqisoper0}
\end{equation}

Formulae (\ref{eqisoper0}) is an immediate consequence of the
definitions. We shall call (\ref{eqisoper0}) the {\em isoperimetric
inequality}. The reader may use the conclusion of this lemma as a
definition of $\kappa _k (S)$. The following upper bound follows by
the inequality $|\partial (\{0\})|\ge \kappa _1$:
\begin{equation}\label{bound}
\kappa _1(S)\le |S|-1.
\end{equation}

The basic intersection theorem is the following:

\begin{theorem}\cite{halgebra,hiso2007}  Assume $|G|\ge 2k-1.$
 Let $A$ be a $k$--atom and let
   $F$   be a   $k$-fragment such that  $|A\cap F|\ge k$. Then
  $A\subset F.$
In particular  distinct $k$-atoms intersect in  at most $k-1$
elements.

\label{inter2frag} \end{theorem}

The structure of $1$--atoms is the following:

\begin{proposition} \label{Cay}\cite{hejc2,hjct}

Let  $ S$  be a  generating subset of an abelian  group $G$ with
$0\in S$.  Let $H$ be a $1$--atom of $S$ with $0\in H$.
   Then
   $H$ is a subgroup. %  generated by $S\cap H$.
 %Moreover for every $1$-fragment $F$ of $S$, $FH=F$.
 Moreover
 \begin{equation}\label{olson}
\kappa _1(S)\geq \frac{|S|}{2},
\end{equation}
%and $ \kappa _1(S)=\frac{|S|}{2}$ holds  if and only if
 %there a subgroup $H$ and an $u$  such that $ S=H\cup H+u$.

\end{proposition}

\begin{proof}
Take $x\in H$. Since $x\in (H+x)\cap H$ and since $H+x$ is a
$1$--atom, we have $H+x=H$ by Theorem \ref{inter2frag}. Therefore
$H$ is a subgroup. Since $S$ generates $G$, we have $|H+S|\ge 2|H|$,
and hence  $\kappa _1(S)=|H+S|-|H|\ge \frac{|S+H|}{2}\ge
\frac{|S|}{2}.$
\end{proof}

Let us mention the following relation between $1$--fragments and
$2$--fragments. We note that a similar relation holds for non
abelian groups and even for abstract graphs.

\begin{lemma}  Let  $S$  be a
finite generating
 $2$--separable subset of an abelian group $G$ with $0\in S$ and $\kappa _2 (S) \leq
 |S|-1$.
 Then $\kappa _2=\kappa _1$. Moreover every $2$--fragment is a
 $1$--fragment. Also every $1$--fragment $F$, with $2\le |F|\le
 |G|-|S|-1$ is a $2$--fragment.
 \label{1f2}
\end{lemma}

Lemma \ref{1f2} follows immediately by the definitions.

 The next result is proved in \cite{Hejcvosp1}. The finite case is reported with almost the same proof in
 \cite{hactaa}.

\begin{theorem} { \cite{Hejcvosp1,hactaa}\label{2atom} Let  $S$  be a finite generating
 $2$--separable subset of an abelian group $G$ with $0\in S$ and $\kappa _2 (S) \leq |S|-1$.
Let
 $H$ be a $2$--atom with $0\in H$. Then $H$ is  a subgroup  or $|H|=2.$

\label{2atomejc} }
\end{theorem}

The next result is an immediate consequence of [\cite{Hejcvosp1},
Theorem 4.6]\label{ejcf}. The finite case ( used to solve Lewin's
conjectures on the  Frobenius number) is reported with almost the
same proof in
 \cite{hactaa}.

\begin{corollary} { [\cite{Hejcvosp1},Theorem 4.6]\label{ejcf}
Let  $S$  be a $2$--separable finite
subset of an
 abelian group $G$ such that $0\in S$, $|S|\leq (|G|+1)/2$ and $\kappa _2 (S) \leq |S|-1$.

 If  $S$ is not an arithmetic progression then there
 is a subgroup $H$ which is a $2$--fragment of $S$.

\label{vosper}                      }
\end{corollary}

\begin{proof}

Suppose that $S$ is not an arithmetic progression.

 Let $H$ be a  $2$- atom such that $0\in H$. If $\kappa
_2\leq |S|-2$, then by Lemma \ref{1f2} $\kappa _2=\kappa _1$ and $H$
is also a $1$--atom. By Proposition \ref{Cay}, $H$ is a subgroup.
Then we may assume $$\kappa _2(S)=|S|-1.$$  By Theorem
\ref{2atomejc}, it would be enough to consider the case $|H|=2$, say
$H=\{0,x\}$. Put $N=\subgp{x}.$

Decompose $S=S_0\cup \cdots \cup S_j$ modulo $N$, where $|S_0+H|\le
|S_1+H| \le \cdots \le |S_j+H|.$ We have $|S|+1=|S+H|=\sum
\limits_{0\le i \le j}|S_i+\{0,x\}|.$

Then $|S_i|=|N|$, for all $i\ge 1$.  We have $j\ge 1$, since
otherwise $S$ would be an arithmetic progression. In particular $N$
is finite.
 We have
$|N+S|<|G|$, since otherwise   $|S|\ge |G|-|N|+1\ge
\frac{|G|+2}{2},$ a contradiction.

Now  \begin{eqnarray*} |N|+|S|-1&=&|N|+\kappa _2(S)\\&\le& |S+N|=
(j+1)|N|\\&\le& |S|+|N|-1, \end{eqnarray*}

 and hence $N$ is a $2$-fragment.
\end{proof}

Corollary \ref{vosper} coincides with [\cite{Hejcvosp1},Theorem
4.6]. A special case of this result is Theorem 6.6 of \cite{hactaa}.
As mentioned in \cite{hplagne1}, there was a misprint in this last
statement. Indeed $|H| + |B| - 1$ should be replaced by $|H| + |B|$
in case (iii) of [ Theorem 6.6, \cite{hactaa}].

The proof of Corollary \ref{vosper} given here uses Proposition
\ref{Cay} and Theorem \ref{2atomejc}. These two results are not
difficult and are proved in around  4 pages [with some possible
simplifications if one forgets about very general results dealing
with  non abelian groups and abstract graphs] in \cite{hiso2007}.

Alternative proofs of Corollary \ref{vosper} (with $|S|\leq |G|/2$
replacing $|S|\leq (|G|+1)/2$),  using Kermperman's Theory were
obtained by Grynkiewicz in \cite{davdecomp} and Lev in
\cite{levkemp}. In the present paper Corollary \ref{vosper} will be
one of pieces leading to a new proof of Kemperman's Theory.

\section{Hyper-atoms }

This section contains the new notion of a hyper-atom. Theorem
\ref{hyperatom} is one of the main results of this paper. As we
shall see later it encodes most of the known results about the
critical pair Theory.
\subsection{Vosper subsets}

Let $0\in S$ be a generating subset of an abelian group $G.$ We
shall say that $S$ is a {\em Vosper subset} if for all $X\subset G$
with $|X|\ge 2$, we have $|X+S|\ge \min (|G|-1,|X|+|S|)$.

Notice that $S$ is a Vosper subset if and only if  $S$ is non
$2$--separable or if $\kappa _2(S)\ge |S|$.

\begin{lemma} { Let $S$ be a finite generating Vosper subset of an
abelian group $G$ such that $0 \in S$.   Let $X\subset G$  be such
that $|X|\ge |S|\ge 3$ and  $|X+S|=|X|+|S|-1$. Then for every $y\in
S$, we have $|X+(S\setminus \{y\})|\ge |X|+|S|-2$. \label{vominus}}
\end{lemma}
\begin{proof}

By  the definition of a Vosper subset. We have $|X+S|\ge |G|-1$.
Then one of the two possibilities:

{\bf Case} 1. $|X+S|= |G|-1$.

 Suppose that $|X+ (S\setminus \{y\})|\le |X|+|S|-3$
and take an element  $z$ of $(X+S)\setminus (X+(S\setminus \{y\}))$.
We have $z-y\in X$. Also $(X\setminus \{z-y\})+S\subset
((X+S)\setminus \{z\})$. In particular we have by the definition of
a Vosper subset, $|(X\setminus \{z-y\})+S|\ge \min
(|G|-1,|X|-1+|S|)=|X|+|S|-1$. Clearly $X+S\supset ((X\setminus
\{z-y\})+S)\cup \{z\}$. Hence $|X+S|\ge |X|+|S|$, a contradiction.

{\bf Case} 2. $|X+S|= |G|$.

 Suppose that $|X+ (S\setminus \{y\})|\le |X|+|S|-3$
and take a $2$--subset   $R$ of $(X+S)\setminus ( X+(S\setminus
\{y\}))$. We have $R-y\subset X$. Also $(X\setminus (R-y))+S\subset
(X+S)\setminus R$. In particular we have by the definition of a
Vosper subset, $|(X\setminus (R-y))+S|\ge \min (|G|-1,|X|-2+|S|)$.
We have  $|X|=1$. Otherwise and since  $X+S\supset ((X\setminus
(R-y))+S)\cup R$, we have  $|X+S|\ge |X|+|S|$,  a contradiction.
Then $|X|=1$. This forces that $|X|=|S|=3$, and hence $|G|=5$. Now
by the Cauchy Davenport Theorem, $|X+(S\setminus \{y\})|\ge
|X|+|S|-2$, a contradiction.
\end{proof}

\subsection{Fragments in quotient groups}

\begin{lemma} { Let $G $ be an abelian group and let $S$ be a
    finite  $2$-separable generating subset containing $0$. Let $H$ be a
subgroup which is a $2$--fragment and let $\phi : G\mapsto G/H$ be
the canonical morphism. Then
\begin{equation}\label{cosetgraph}
\kappa _1(\phi (S))=  |\phi (S)|-1.
\end{equation}
 Let $K$ be a subgroup which is a $1$--fragment of
$\phi (S)$. Then $\phi ^{-1}(K)$ is a $2$--fragment of $S$.
}\label{quotient}
\end{lemma}
\begin{proof}

 Put $|\phi (S)|=u+1$.  Since $|G|>|H+S|,$ we have $\phi (S)\ne G/H$,
 and hence $\phi (S)$ is
$1$--separable.

Let $X\subset G/H,$ be such that  $X+\phi (S)\neq G/H$. Clearly
$\phi^{-1} (X)+S\neq G$. Then $|\phi^{-1} (X)+S|\ge |\phi^{-1}
(X)|+\kappa _1(S)= |\phi^{-1} (X)|+u|H|.$

It follows that $|X+\phi (S)||H|\ge |X||H|+u|H|.$ Hence $\kappa
_1(\phi (S))\ge u=|\phi (S)|-1$. The reverse inequality is obvious
and follows by (\ref{bound}). This proves (\ref{cosetgraph}).

Let $K$ be a subgroup which is a $1$--fragment of $\phi (S)$. Then
$|K+\phi (S)|=|K|+u$. Then $|\phi ^{-1} (K)+S|=|K||H|+u|H|.$ In
particular $|\phi ^{-1} (K)$ is a $2$--fragment.\end{proof}

\subsection{The fundamental property of hyper-atoms}

Let $S$ be a finite generating subset of an abelian group $G$ such
that $0 \in S.$ Theorem \ref{Cay} states that there is a $1$--atom
of $S$ which is a subgroup. A subgroup with maximal cardinality
which is a $1$--fragment will be called a {\em hyper-atom}. This
definition may adapted to non-abelian groups and even abstract
graphs. As we shall see the hyper-atom is more closely related to
the critical pair theory than the $2$--atom.

\begin{theorem}\label{hyperatom}
Let $S$ be a finite generating subset of an abelian group $G$ such
that $0 \in S,$  $| S | \leq (|G|+1)/2$ and $\kappa _2 (S)\le
|S|-1.$ Let $H$ be a hyper-atom of $S$. Then

\begin{itemize}
  \item [(i)] $\phi (S)$ is either an arithmetic progression or a Vosper
subset, where $\phi$ is the canonical morphism from $G$ onto
$G/H$.
 \item [(ii)] Let $X\subset G/H$  be such that $|X+\phi (S)|=|X|+|\phi (S)|-1$.
Then for every $y\in \phi (S)$, $|X+(\phi (S)\setminus y)|\ge
|X|+|\phi (S)|-2$.
\end{itemize}
\end{theorem}

\begin{proof}

Let us show that $2|\phi (S)|-1\le \frac{|G|}{|H|}.$ Clearly we may
assume that $G$ is finite.

Observe that $2|S+H|-2|H|\le 2|S|-2< |G|.$ It follows, since$|S+H|$
is a multiple of $|H|$, that $2|S+H|\le  |G|+|H|,$ and hence $2|\phi
(S)|\le \frac{|G|}{|H|}+1.$

 Suppose now that  $\phi (S)$ is not a Vosper subset. By the
 definitions $\phi (S)$ is $2$--separable and $\kappa _2(\phi(S))\le |\phi(S)|-1.$

 Observe that  $\phi(S)$ can not
have  a $2$-fragment $M$ which is a subgroup. Otherwise  by Lemmas
\ref{quotient} and  \ref{1f2}, $\phi ^{-1}(M)$ is a $2$--fragment of
$S$ containing strictly $H$, contradicting the maximality of $H$. By
Corollary \ref{vosper}, $\phi(S)$ is an arithmetic progression.

Now (ii) holds by Lemma \ref{vominus} if $\phi (S)$ is a Vosper
subset. It is also obvious if $\phi (S)$ is an arithmetic
progression.
\end{proof}

\begin{corollary}\label{plagne}

Let S be a generating subset of a finite abelian group G such that
$0\in S$ and  $|S| \le \frac{|G|}2$, then one of the following
conditions holds:
\begin{itemize}
  \item [(i)] $S$ is an arithmetic progression,
    \item [(ii)] there is a subgroup $H\neq \{0\}$
 such that $|H + A| < \min(|G| - 1, |H| + |S|)$,
   \item [(iii)] for any $X$
 such that $|X | \ge 2$, $|S + X| \ge  \min(|G| - 1, |S| + |X
|)$.
\end{itemize}
\end{corollary}
Notice that the main aim of the authors of Corollary \ref{plagne}
was to give an application to sum free sets in finite abelian
groups. The infinite case was irrelevant for this purpose. However
the proof works if $S$ is a   finite subset of an  abelian group if
one uses [\cite{Hejcvosp1}, Theorem 4.6] instead of   Theorem 6.6 of
\cite{hactaa}. Alternative proofs of Corollary \ref{plagne} using
Kermperman's Theory were obtained  by Grynkiewicz in
\cite{davdecomp} and Lev in \cite{levkemp}.

Theorem \ref{hyperatom} implies clearly Corollary \ref{plagne} with
some improvements:
\begin{itemize}
\item
The subgroup $H$ in Theorem \ref{hyperatom}  is well described
as a hyper-atom;
\item
We have also an equality $|H+S|-|H|=\kappa _1$,  much precise
than the inequality $|H+S|\le |H|+|S|-1$.  This equality will be
needed later;
\item The condition $|S| \le \frac{|G|}2$ is relaxed to $|S| \le
\frac{|G|+1}2$.
\end{itemize}

Part (ii) of Theorem \ref{hyperatom} is  a critical pair result of a
new type, that will be used later to prove Kemperman's structure
Theorem.

\section{Quasi-periodic decompositions}

\begin{theorem}\label{quasiperiod}
Let $S,T$ be   finite  subsets  of
 an abelian group $G$ with $|S+T|=|S|+|T|-1$.

Assume moreover that   $S+T$ is aperiodic. Then one of the following
holds:
\begin{itemize}
\item [(i)]  $S$ and $T$ are $K$-quasi-periodic, for some nonzero subgroup $K$.
\item [(ii)]  The pair $\{S,T\}$ is a strict elementary pair.

\end{itemize}

\end{theorem}

\begin{proof}

The proof is by induction on $|S|+|T|$, the result being obvious for
$|S|+|T|$ small.
 We may assume  clearly that $0\in S$. We may assume $\min (|S|,|T|)\ge 2$, since otherwise $\{S,T\}$
is a strict elementary pair and (ii) holds. Without loss of
generality we may assume $2\le |S|\le |T|.$

{\bf Claim 1} If $T\not\subset \subgp{S}$, then the result holds.

\begin{proof}
Decompose $T=\bigcup _{i\in U} T_i$ modulo  $\subgp{S}$. By
(\ref{olson}), $\kappa _1(S)\ge \frac{|S|}2.$ Put $V=\{i\in U :
|T_i+S|<|\subgp{S}|\}.$ By (\ref{eqisoper0}) we have

\begin{eqnarray}
|T+S|
&\ge&(|U|-|V|)|\subgp{S}|+ \sum \limits_{i\in V}
|T_i+S|\label{reduct}\\
 &\ge&(|U|-|V|)|\subgp{S}|+ \sum \limits_{i\in V}
|T_i|+|V|\frac{|S|}2\ge |T|+|V|\frac{|S|}2.\nonumber
\end{eqnarray}

It follows that $|V|\le 1.$ But $|V|\ge 1$, since otherwise
$T+S=T+S+\subgp{S}$ . Put $V=\{\omega\}$. By Kneser's Theorem
$|T_{\omega}+S|\ge |T_{\omega}|+|S|-1.$ By (\ref{reduct}) we have
\begin{eqnarray*} |T|+|S|-1=|T+S| &\ge&
(|U|-1)|\subgp{S}|+|T_{\omega}|+|S|-1
\end{eqnarray*}
Therefore Then $S$ and $T$ are $\subgp{S}$-quasi-periodic.
\end{proof}

By Claim 1, we may assume without loss
  of generality that
  $$G=\subgp{S}.$$
 We may assume that $S$ is not an arithmetic
 progression since otherwise $T$ would be an arithmetic progression
 with the same difference, and (ii) would be satisfied.

Assume first $|G|-|T+S|=|T^S|< |T|$. Then $G$ is finite. Observe
that $T^S-S$ is aperiodic, otherwise by Lemma $T+S=(G\setminus
({T^S}^{-S}))+S$ would be periodic. By Kneser's Theorem $|T^S-S|=
|T^S|+|S|-1$. By the definition $(T^S-S)\cap T=\emptyset$. Therefore
$|T^S-S|\le |G|-|T|=|G|-|S+T|+|S+T|-|T|\le |T^S|+|S|-1$. Hence
$T^S-S=G\setminus T$, and hence ${T^S}^{-S}=T$.

 Then one of the
following conditions holds by the induction hypothesis:

\begin{itemize}
  \item
 $S,T^S$ are $N$--quasi-periodic, for some non zero subgroup
$N$. Therefore $T=G\setminus (T^S-S)$ is $N$--quasi-periodic.
The result holds in this case.

 \item   The pair $\{S,T^S\}$ is an elementary pair. Observe that  $S$ is not an arithmetic progression
 and hence (SP1) can not be satisfied for the pair $\{S,T^S\}$.
 Also $|T^S-S|=|G|-|T|\ge 2$ and then (SP3) can not be satisfied
 for the pair $\{S,T^S\}$.

  Then necessarily is $|T^S|=\min (|T^S|,|S |)=1.$ Let $c$
 denotes the unique element of $G\setminus (T+S)$. Then
$c-T\subset G\setminus S$. But $|c-T|=|T|\ge
|S+T|-|S|+1=|G|-|S|.$ This shows that $c-T= G\setminus S.$
Observe that $T$ is aperiodic, since otherwise $T+S$ would be
periodic.

 {\bf Case 1}: $|(c-S)\cap T|\neq 1$ for every $c\in G$. In this
case $\{T,S\}$ is a strict elementary pair and (ii) holds.

{\bf Case 2}: $|(c-S)\cap T|= 1$ for some $c\in G$. Put
 $c=x_1+y_1$, where $x_1\in T$ and $y_1\in S$. Put
$T'=T\setminus \{x_1\}$. Clearly $|T'+S|\le |H|-2=|T'|+|S|-1$.
Let $Q$ denotes the period of $T'+S$.  By Kneser's Theorem
\ref{kneser}, $|H|-2\ge |T'+S|\ge |T'+Q|+|S+Q|-|Q|$. This forces
that $|Q|=1$, since otherwise (observing that $S$ is aperiodic)
we have  $|T|+|S|-1-|Q|=|H|-|Q|\ge |T'+S|\ge |T'+Q|+|S+Q|-|Q|\ge
(|T|-1)+(|S|+1)-|Q|$, a contradiction.

Therefore  $|H|-2\ge |T'+S|\ge |T'|+|S|-1=|H|-2.$  Put
$\{x_1,x_2\}=H\setminus (T'+S)$ and $d=x_2-x_1$. Since
$x_2-d=x_1$, we have $x_2\notin T'+S+d$. Since $T'+S$ is
aperiodic, we have  $|H|-1=|T'+S|+1\le |T'+S+\{0,d\}|\le |H|-1$.
 Since $T'+S+\{0,d\}$ is aperiodic, we have Kneser's Theorem
\ref{kneser} $|(S+\{0,d\})+T'|\ge |S+\{0,d\}|+|T|-2.$ It follows
that $|S+\{0,d\}|\le |S|+1$. Hence $S$ is
$\subgp{d}$--quasi-periodic. Similarly $T$ is
$\subgp{d}$--quasi-periodic is $\subgp{d}$--quasi-periodic. Then
(i) holds in this case.

\end{itemize}
 So we may assume that
$|T|\le |T^S|.$

By our assumptions $|T^S|=|G|-|T+S|\ge |T|\ge |S|$, we have

 \begin{eqnarray*}
 3|S+T|&= &2|S+T|+|S|+|T|-1\\
 &\le&|G|-|S|+|G|-|T|+|S|+|T|-1= 2|G|-1,
  \end{eqnarray*}
In particular we have

\begin{equation}\label{eq2n/3}
|S+T|\le   \frac{2|G|-1}{3}.
\end{equation}

Let $H$ be a hyper-atom of $S$ and let $\phi :G\mapsto G/H$ denotes
the canonical morphism. Put $|\phi(S)|=u+1$ and $|\phi(T)|=t+1$. Put
$q=\frac{|G|}{|H|}$.

Take a $H$--decomposition $S=\bigcup \limits_{0\le i\le u}S_i$ such
that $|S_0|\ge \cdots \ge |S_u|$. By the definition we have
$u|H|=|H+S|-|H|=\kappa _1\le |S|-1.$ It follows that for all $u\ge
j\geq 0$
\begin{equation}\label{plein}
|S_{u-j}|+\cdots +|S_u|\ge j|H|+1
\end{equation}

It follows that $|S_0| \ge \frac{|H|+1}{2}$. In particular $S_0$
generates $H$.
% More generally $S_i-S_i$ generates $H$ for all $0\le i \le u-1$.
We shall use this fact in the application of the isoperimetric
inequality.

 Take a $H$--decomposition $T=\bigcup \limits_{0\le i\le t}T_i$.

 By (\ref{cosetgraph}), $\kappa _1
(\phi(S))=|\phi (S)|-1=u.$ Put $\ell =\min (q-t-1,u)$.

 By Proposition \ref{strongip} applied to  $\phi (S)$ and $\phi (T)$, there is a subset $J\subset
[0,t]$ with cardinality $\ell$ and a family $\{ mi ;i\in J\}$ of
integers in $[1,u]$
  such that  $T+S $ contains the $H$--decomposition $(\bigcup \limits_{0\le i\le
  t}T_i+S_0)\cup (\bigcup \limits_{ i\in
  J}T_i+S_{mi})$.

Put $R=(S+T)\setminus  ((\bigcup _{i\in J} {T_i+S_{mi}}+H)\cup
(\bigcup _{0\le i\le t} T_i+H))$.

We shall choose such a $J$ in order to maximize $|J\cap P|.$ We
shall write $E_i=(S+T)\cap (T_i+H)$, for every $i\in [0,t]$. Also we
write $E_{mi}=(S+T)\cap (T_i+S_{mi}+H)$, for every $i\in J$.

 We put also  $W=\{i \in [0,t] :  |E_i|<|H|\},$ and $P=[0,t]\setminus W.$

Since $|T|\ge |S|$ we have $|T+H|\ge |S|>\kappa_2(S)=u|H|.$ It
follows that $t+1=|\phi (T)|\ge u+1.$ Then $t+1-|J|>0.$ In
particular $I\neq \emptyset ,$ where $I=[0,t]\setminus J.$

Let $X$ be a subset of $I$ and  let $Y$ be a subset of $J$. We have
\begin{eqnarray}
|S+T|-|R|
&\ge &\sum \limits_{i\in X\cup Y}|E_i|+ \sum \limits_{i\in I\setminus X\cup J\setminus Y}|T_i+S_{0}|
+\sum\limits _{i\in J\setminus Y}|T_i+S_{mi}|+\sum \limits_{i\in Y}|E_{mi}|\nonumber\\
&\ge &\sum \limits_{i\in X\cup Y}|E_i|+ \sum \limits_{i\in I\setminus X \cup J\setminus Y} |T_i| + (u-|Y|)|S_{0}|+\sum \limits_{i\in Y}|E_{mi}|\label{debut}\\
&\ge &\sum \limits_{i\in X\cup Y}|E_i|+ \sum \limits_{i\in I\setminus X\cup J\setminus Y} |T_i|+  (u-|Y|)|S_{0}|
+|Y||S_{u}|\label{eqdebut}
\end{eqnarray}

Put $F=\{i\in I\cap P : (T_i+S) \cap (\bigcup _{i\in W}T_i+H)\neq
 \emptyset\}$.

We shall use the following obvious facts: For all $i\in W$, we have
by (\ref{olson}), $|E_i|\ge|T_i+S_0|\ge |T_i|+\kappa _1(S_0)\ge
|T_i|+\frac{|S_0|}{2}.$ For every $i\in F$, $T_i+S_{ri} \subset
T_j+H$ for some $1\le ri\le u$ and some $j\in W.$ Hence we have
$|T_i|+|S_{u}|\le |T_i|+|S_{ri}|\le |H|=|E_i|,$ by Lemma
\ref{prehistorical}.

Let $U$ be a subset of $W\cap J$. Put $X=I$ and $Y=U$. By
(\ref{eqdebut}), we have

\begin{eqnarray}
|S+T|-|R|
&\ge &\sum \limits_{i\in U\cup (W\cap I)\cup (P\cap I)}|E_i|+ \sum \limits_{i\in J\setminus U} |T_i|+ (u-|U|)|S_{0}|+|U||S_{u}|
\\
&\ge &\sum \limits_{i\in (P\cap I)\setminus F}|T_i|+\sum \limits_{i\in F}(|T_i|+|S_u|)+ \sum \limits_{i\in (W\cap I) \cup U}(|T_i|+\frac{|S_{0}|}2)
+|J\setminus U||S_{0}|
+|U||S_{u}|\nonumber\\
&\ge& |T|+  |J\setminus U||S_0|+(|U|+|F|)|S_u|+|(W\cap I) \cup U|\frac{|S_{0}|}2.\label{eqvide}
\end{eqnarray}

{\bf Claim 2} $q \ge |\phi (S)|+|\phi (T)|-1$, and hence $\ell =u$.

\begin{proof} The proof is by contradiction. Suppose that $q < |\phi (S)|+|\phi
(T)|-1$.

 Assume first $u\ge 2$. By Lemma \ref{prehistorical}, the
are two distinct values of the pair $(s,t)$ such that $T_s+S_t
\subset E_{mi}$, for every $i\in J$. In particular $|E_{mi}|\ge
|S_{u-1}|$, for every $i\in J$. Also $|E_{i}|\ge |S_{0}|$, for every
$i\in [0,t]$.

  Observe that $2t> t+u\ge q$. We have using (\ref{plein})

  $ 2|S_0|\ge |S_0|+|S_{u-1}|\ge \frac{2}{3}(|S_u|+|S_{u-1}|+|S_{u-2})>\frac{4|H|}3$.
  By (\ref{eqdebut}), applied with $X=I$ and $Y=J$, we have

\begin{eqnarray*}
|S+T|&\ge& \sum \limits_{0\le i \le t} |S_0|+ \sum \limits_{i\in J} |S_{u-1}|= (t+1)|S_0|+(q-t-1)|S_{u-1}|\\
&=& (2t+2-q)|S_0|+(q-t-1)(|S_0|+|S_{u-1}|)\\
&>& (2t+2-q)\frac{2|H|}3+\frac{4|H|(q-t-1)}3=\frac{2|G|}3,
\end{eqnarray*}

contradicting (\ref{eq2n/3}).

Assume now $u=1.$ From the inequality $|T+S|\le |T|+|S|-1$, we see
that $\kappa _1(S)\le |S|-1$. Therefore we have by (\ref{eq2n/3}),
$\frac{2|G|}{3}>|T+S|\ge |T|+\kappa _1(S)\ge |S|+|H|>2|H|$, and
hence
$$q\ge 4.$$

 We have $(t+1)+(u+1)-1<|\phi (S+T)|\le  q.$ Then $t+1=q.$
Hence $\ell =|J|=0$.
 We have $|W|\geq 1,$ since
otherwise $G= T+H\subset S+T$. We have $|W|\leq 3,$ by
(\ref{eqvide}) applied with $U=\emptyset$. Therefore $|P|\geq
t+1-3\ge 4-3=1$. There is clearly   $i\in P$ with $T_i+S_1 \subset
T_j+H$ for some $j\in W,$ and hence $|F|\ge 1$. By (\ref{eqvide})
applied with $U=\emptyset$, $|T+S|\ge |T|+|W|\frac{|S_0|}{2}+|S_1|$,
and hence $|W|\le 1$. It follows that $|S+T|\ge |G|-|H|=
|G|-\frac{|G|}q \ge \frac{3|G|}{4}$, contradicting
(\ref{eq2n/3}).\end{proof}

We must have $R = \emptyset$, since otherwise by (\ref{eqvide})
applied with $U=\emptyset$, $|S+T|-|R|\ge |S+T|-|S_u||\phi(R)|\ge
|T|+u|S_0|+|S_u|\ge |T|+|S|,$ a contradiction. In particular

\begin{equation} \label{cdcp}
|\phi (S+T)| = |\phi (S)|+|\phi (T)|-1.
\end{equation}

 {\bf Claim 3}.  $J \cap P \neq \emptyset.$

\begin{proof}
Suppose the contrary and take $k\in J\cap W$. Put $U=\{k\}$. By
(\ref{eqvide}),

\begin{eqnarray*}
|S|+|T|>|S+T|&\ge& |T|+ (u-1)|S_{0}|+|S_{u}|+(|W\cap I|+1)\frac{|S_{0}|}2.
\end{eqnarray*}

It follows that $I\subset P$. Since $S$ generates $G$, we have
$|\bigcup _{i\in I} T_i+H+S|>|\bigcup _{i\in I} T_i+H|$.

 We must have $ (\bigcup _{i\in I} T_i+H+S)\cap (\bigcup _{i\in J} E_{mi}+H)=\emptyset $,
 since otherwise by replacing a suitable element of $J$ with some $p\in I$, we
may increase strictly $ |J\cap P|,$ observing that $I\subset P$.

By (\ref{cdcp}), there are $i\in I$, $j\in J$ and $p\in [1,u]$ such
that $ T_i +S_p$ is congruent $T_j+S_{mj}$. It follows that $F\neq
\emptyset$.

By (\ref{eqvide}) applied with $U=\emptyset$,

\begin{eqnarray*}
|S+T|&\ge& |T|+ u|S_{0}|+|S_{u}|\ge |T|+|S|,
\end{eqnarray*}
a contradiction proving the claim.
\end{proof}

Take $r\in  J$ with $|E_r|=|H|$. Such an $r$ exists by Claim 3.

 {\bf Claim 4} $T_i+H+S_j=T_i+S_j,$
for all $0\le j\le u-1$.

 \begin{proof} By Lemma \ref{prehistorical}, it would be enough to show the
following:
\begin{equation} \label{FINAL}
|T_k|+|S_{u-1}|>|H|,
\end{equation}
for every $k\in [0,t]$.
 Suppose the contrary.

Notice that $|E_{mr}|\ge \max (|T_r|,|S_u|)$ and that $|E_k|\ge
|S_0|$. Also $|T_k|+|S_{u-1}|\le |H|=|E_{mr}|$ by our choose of $r$.
We shall use these inequalities and (\ref{debut}) with
$X=\{k,r\}\cap I$ and $Y=\{k,r\}\cap J$.

 By (\ref{debut}) we have for $k\neq r$,
\begin{eqnarray*} |S+T|
&\ge&  |T|-|T_k|-|T_r|+ (u-|X|)|S_0|+|T_k|+|S_{u-1}|+|S_{0}|+|T_r|+|Y||S_{u}|\\&\ge&  |T|+(u-1) |S_0|+|S_{u-1}|+|S_{u}|\ge |T|+|S|,
\end{eqnarray*}
leading a contradiction. If $k=r$ the contradiction comes more
easily.
\end{proof}

Since $|S+T|< |G|$, we must have by Lemma  \ref{prehistorical},
$$2|S|\le |S|+|T|\le |G|.$$
Now by (\ref{cdcp}) and Theorem \ref{hyperatom},
$|\phi(T+(S\setminus S_u))|\ge t+u$. Take a subset $\Omega $ of
$\phi(T+(S\setminus S_u))$ with $|\Omega |=u+t$. By (\ref{cdcp}), $
\phi(T+S) =\Omega \cup \{\omega\}$, for some $\omega \in G/H$. By
Claim 4, $(\phi ^{-1}(\Omega))\cap (T+S)$ is  $H$--periodic.
Necessarily there is $s$ such that $|E_{ms}|<|H|.$ Then by Claim 4
$E_{ms}=T_s+S_u$. Since $T+S$ is aperiodic, and since
$(T+S)\setminus E_{ms}$ is $H$--periodic, we have that $T_s+S_u$ is
aperiodic. By Kneser's Theorem, $|T_s+S_u|\ge |T_s|+|S_u|-1.$ Now we
have
\begin{eqnarray*} |S|+|T|-1 &=&|S+T|\\ &=& |(\psi ^{-1}(\Omega))\cap
(T+S)|+|E_{ms}|\\
&\ge&  (t+u)|H|+|E_{ms}|\\&\ge& (t+u)|H|+|T_{s}|+|S_{u}|-1\ge |T|+|S|-1.
\end{eqnarray*}

Therefore $|T|=t|H|+|T_{s}|$ and $|S|=u|H|+|S_{u}|$. Hence $T$ and
$S$ are $H$--periodic.\end{proof}

\begin{proofof}{Theorem}{\ref{Kemperman}}

The implication (II) $\Rightarrow $ (I) is quite easy. Let us prove
the implication (I)$\Rightarrow $(II). Suppose that (I) holds.

%"Yo Allah yobe modyore"
Assume first that  $A+B$ is aperiodic. Note that $(\emptyset , A)$
and $(\emptyset , B)$ are $G$-quasi-periodic decompositions. Take a
subgroup $H$ with minimal cardinality $|H|\ge 2$ for which there are
$H$--quasi-periodic decompositions $A=A_0\cup A_1$ and $B=B_0\cup
B_1$. Let $\phi : G\mapsto G/H$  be the canonical morphism. Take
$a_1\in A_1$ and $b_1\in B_1$.

Notice that $A_1$ and $B_1$ have no $P$-quasi periods for some $2\le
|P|<|H|$, otherwise $|H|$ would be not minimal. By Theorem
\ref{quasiperiod}, the pair $\{A_1,B_1\}$ is an elementary pair.
Since $A+B$ is aperiodic, $\phi(a_1)+\phi (b_1)$ has a unique
expression.

Assume now that  $A+B$ is periodic.

 Let $H$ be a  period of
$A+B$ with a prime order and let $\phi : G\mapsto G/H$  is the
canonical morphism.

Let $C$ denotes the set of elements of $A+B$ having a unique
expression. Clearly $c\in C$. To each $x\in C$, choose $a_x\in A$
and $b_x\in B$ such that $x=a_x+b_x$. Put $A_x=A\cap (a_x+H)$.

  Observe that $\phi (c)=\phi (a_c)+\phi (b_c)$ has a unique
expression. Hence by Scherck's Theorem \ref{scherk}, $|\phi(A)+\phi
(B)|\ge |\phi(A)|+|\phi (B)|-1.$ We must have $|\phi(A)+\phi (B)|=
|\phi(A)|+|\phi (B)|-1,$ since otherwise
$|A+B|=|\phi(A+B)||H|=|A|+|B|.$
 By
Lemma \ref{prehistorical}, we have $$|A_x|+|B_x|\le |H|+1.$$ Observe
that $|A|+|B|-1 = |A+B| = |\phi(A+B)||H|=|A+H|+|B+H|-|H|.$ It
follows that the trace of $A$ (resp.  $B$)on any coset $\neq A_x+H$
(resp.   $\neq B_x+H$) has cardinality $=|H|$. It follows that
$A\setminus A_x$ and $B\setminus B_x$ are $H$-periodic sets.

If
$|C|=1$, then $\{A_c,B_c\}$ is an elementary pair. So we may assume
that $|C|\ge 2.$ We may assume that  $|A_c|\ge 2$, since otherwise
$\{A_c,B_c\}$ is an elementary pair (verifying SP2). Assume first
that $|A_c|= 2$, say $A_c=\{a,a+d\}$. Then $|B_c|=|H|-1$ and hence
$B_c$ is an arithmetic progression of difference $d$ (recall that
$|H|$ is a prime). It follows that $\{A_c,B_c\}$ is an elementary
pair (verifying SP1). So we may assume that $|A_c|\ge 3$.

Suppose that there is $v\in C\setminus \{c\}$, with $A_v\neq A_c$ or
$B_v\neq B_c$. Without loss of generality we may take $A_v\neq A_c$.
Since $A\setminus A_v$ is $H$--periodic, we have  $|A_c|=|H|$. Then
$|B_c|=|H|+1-|A_c|=1$. It follows that $\{A_c,B_c\}$ is an
elementary pair.

Therefore there is $v\in C\setminus \{c\}$, with $A_v=A_c$ and
$B_v=B_c$. Put $c=c_1+d_1$ and $v=c_2+d_2$, with $c_1,c_2\in A$ and
 $d_1,d_2\in B.$
 Now we have $|B_c|\ge |A_c|\ge 3.$ If $|H|=5$, then $A_c, B_c$,
have cardinality $=3$, and hence they are arithmetic progressions.
The unique expressibility of $c$ implies that $A_c, B_c$ have the
same difference. It follows that $\{A_c,B_c\}$ is an elementary pair
(verifying SP1). So we may assume that $|H|\ge 7$. Therefore
$2|B_c|\ge |B_c|+|A_c|=|H|+1\ge 8.$

 Put $B'=B_c\setminus \{c_1,d_1\}$. We have clearly $|H|-2=|A'|+|B_c|-1\ge |B'+A_c|$.
By Vosper's Theorem  $A_c, B'$ are arithmetic progressions with a
same difference say $d$. Since $|(B_c\setminus \{c_i\})+A_c|\le
|B_c\setminus \{c_i\}|+|A_c|-1$, we must have that $B_c\setminus
\{c_i\}$ is an arithmetic progression with difference $=d$ and
extremity $c_i$. This shows that $\{A_c,B_c\}$ is an elementary pair
(verifying SP1). theorem is proved in this case.

\end{proofof}

%The above lemma is proved implicitly  by Kemperman
%[\cite{kempacta},Lemma 5.2]. Our proof looks easier. Another proof
%follows from Lev's [\cite{levkemp}, Appendix: Propositions 1 and 2].

\section{Appendix: A short proof of Scherck's Theorem }

\begin{proofof}{Theorem}{\ref{scherk}}

We start by a special.

{\bf Claim} For any two finite subset $A,B$ such that $A\cap
(-B)=\{0\}$, we have $|A+B|\ge |A|+|B|-1.$

 The proof is by induction on $|B|$, the
result being obvious for $|B|=1$.

 Since $0\in A\cap B,$ we have
$|A+B|\ge |A\cup B|=|A|+|B|-|A\cap B|$. Therefore we may assume that
there is $b \in A\cap B,$ with $b\neq 0.$

Put $B'=B\cap (B-b)$. We have $B-b\neq B,$ since otherwise we would
have $b\in A\cap -B,$ a contradiction. Hence $0\in B'$ and
$|B'|<|B|$.

Put $A'=A\cup A+b$.  We have $A'\cap -B'= (A\cap -B)\cap (-B+b) \cup
A+b\cap (-B+b)\cap (-B)=\{0\}.$

The result follows by induction.

We may put $c=c_1+c_2$, where $c_1\in X$ and $c_2\in Y$. Put $
A=X-c_1$ and $B=Y-c_2$. The result follows by the Claim.
\end{proofof}


\begin{thebibliography}{99}



\bibitem{balart} E. Balandraud, Un nouveau point de vue isop\'erimetrique
appliqu\'e au th\'eor\`eme de Kneser, {\it Preprint}, december 2005.
%\bibitem{balthes} E. Balandraud, Quelques r\'esultats combinatoires en
%Th\'eorie Additive des Nombres, {\it Th\`ese de Doctorat de
%l'Universit\'e de Bordeaux I}, May 2006.




\bibitem{cauchy} A.  Cauchy,  Recherches  sur  les  nombres,  {\it J.  Ecole  polytechnique}
9(1813), 99-116.



\bibitem{davenport} H. Davenport, On the addition of residue  classes, {\it J.  London  Math.
Soc.} 10(1935), 30--32.


%\bibitem{desfreim} J-M. Deshouillers, G. A. Freiman, A step beyond
 % Kneser's Theorem. {\it Proc. London Math. Soc.} (3)86 (2003), no 1,
 % 1--28.



%\bibitem{dixmier}  J. Dixmier, Proof of a conjecture
%by Erd\"os, Graham concerning the problem of Frobenius, {\it J.
%number Theory} 34 (1990), 198-209.


\bibitem{gerlodinger} A. Geroldinger, F. Halter-Koch,  Non-unique factorizations.
Algebraic, combinatorial and analytic theory. Pure and Applied
Mathematics (Boca Raton), 278. Chapman \& Hall/CRC, Boca Raton, FL,
2006. xxii+700 pp.


%\bibitem{greenruz} B. Green, I. Z. Ruzsa,
 %Sets with small sumset and rectification. {\it Bull. London Math. Soc.} 38 (2006), no. 1, 43--52.

\bibitem{davdecomp} D. Grynkiewicz, Quasi-periodic decompositions and the Kemperman's structure
theorem,  European J. Combin.  26  (2005),  no. 5, 559--575.

\bibitem{davkem} D. Grynkiewicz, A step beyond Kemperman's structure Theorem, Preprint May 2006.

%\bibitem{hall} M. Hall, Distinct representatives of subsets, {\it Bull. Amer. Math. Soc.}
%54 (1948), 922--926.

\bibitem{hcras} Y.O. Hamidoune, Sur les atomes d'un graphe orient\'e,
 {\it C.R. Acad. Sc. Paris A}  284 (1977),   1253--1256.


\bibitem{hjct} Y.O. Hamidoune,   Quelques probl\`emes de connexit\'e dans les
graphes orient\'es,  {\it J. Comb. Theory} B 30 (1981), 1-10.


\bibitem{hejc2} Y.O. Hamidoune, On the connectivity of Cayley digraphs,
{\it Europ. J.  Combinatorics}, 5 (1984), 309-312.

\bibitem{Hejcvosp1} Y.O. Hamidoune, Subsets with small sums in abelian groups I: The Vosper property.
European J. Combin. 18 (1997), no. 5, 541--556.


\bibitem{halgebra} Y.O. Hamidoune, An isoperimetric method in additive theory.
{\it J. Algebra} 179 (1996), no. 2, 622--630.


 %\bibitem{hast} Y.O. Hamidoune,  On small subset product in a group.
%Structure Theory of set-addition,  {\it Ast\'erisque}  no. 258(1999),
%xiv-xv, 281--308.

 \bibitem{hactaa} Y.O. Hamidoune, {Some results in Additive number
Theory I: The critical pair Theory}, Acta Arith. 96, no. 2(2000),
97-119.



\bibitem{hiso2007} Y.O. Hamidoune, Some additive applications of the isopermetric approach,
 http://arxiv.org/abs/math./07060635.


\bibitem{hglobal} Y.O. Hamidoune,
  The global isoperimetric methodology applied to Kneser's Theorem, Preprint August-2007.
% http://arxiv.org/abs/math./07060635.

\bibitem{hplagne1} Y. O. Hamidoune , A. Plagne.
{A new critical pair theorem applied to sum-free sets}.
  {\it Comment. Math. Helv.} 79 (2004), no. 1, 183--207.



\bibitem{hrodseth2} Y. O. Hamidoune, {\O}. J. R{\o}dseth, An inverse theorem modulo $p$,
{\it Acta Arithmetica} 92 (2000)251--262.



%\bibitem{hgoaa1}Y. O. Hamidoune, O. Serra, G. Z\'emor,
%{On the critical pair theory in $\zp$}, {\it Acta Arith.} 121 (2006), no. 2, 99--115.


%\bibitem{hgochowla}Y. O. Hamidoune, O. Serra, G. Z\'emor,
%{On the critical pair theory in Abelian groups: Beyond Chowla's Theorem}, Preprint september 2006.\\
%\verb=http://arxiv.org/abs/math.NT/0603478




%\bibitem{jin} R. Jin, Solution to the inverse problem for upper asymptotic density.  {\it J. Reine Angew. Math.}  595  (2006), 121--165.



\bibitem{kempacta} J. H. B. Kemperman, On small sumsets in Abelian groups,
{\it Acta Math.} 103 (1960), 66--88.

%\bibitem{kempcompl} J.H.B. Kemperman, On complexes in a semigroup,  {\it Nederl. Akad. Wetensch. Proc. Ser. A.} 59=  {\it Indag. Math.} 18(1956), 247--254.

%\bibitem{knesrcomp} M. Kneser, Summenmengen in lokalkompakten  abelesche Gruppen,
%{\it Math. Zeit.} 66 (1956), 88--110.


\bibitem{levkemp} V. F. Lev,  Critical pairs in abelian groups and Kemperman's structure theorem.
{\it Int. J. Number Theory} 2 (2006), no. 3, 379--396.

%\bibitem{manlivre} H.B. Mann, {\it Addition Theorems},   R.E.
%Krieger, New York, 1976.
%\bibitem{manams} H. B. Mann,  An addition theorem for sets of elements of an Abelian group,{\it Proc. Amer. Math. Soc.} 4 (1953), 423.
\bibitem{natlivre}M. B. Nathanson,
{\it Additive Number Theory. Inverse problems and the geometry of
sumsets}, Grad. Texts in Math. 165, Springer, 1996.



%\bibitem{olsonjnt} J.E. Olson, On the sum of two sets in a group,
 % {\it J. Number Theory} 18 (1984), 110--120.

\bibitem{scherck}P. Scherk, L.Moser,   Advanced Problems and Solutions: Solutions: 4466,{\it  Amer. Math. Monthly}  62  (1955),  no. 1, 46--47.
%\bibitem{serra} O. Serra,  An isoperimetric method for the small sumset problem. {\it Surveys in combinatorics} 2005, 119--152,
%{\it London Math. Soc. Lecture Note Ser.}, 327, Cambridge Univ. Press, Cambridge, 2005.
%\bibitem{SZ}  O. Serra, G. Z\'emor,
%On a generalization of a theorem by Vosper,{\it Integers} 0, A10,
%(electronic) 2000.\\
%\verb=http://www.integers-ejcnt.org/vo10.html



\bibitem{tv} T. Tao, V.H. Vu,  {\it Additive Combinatorics}, Cambridge Studies
in Advanced Mathematics 105 (2006), Cambridge University Press.




%\bibitem{vosper1} G. Vosper, The critical pairs of subsets of a group of prime order,
 %{\it J. London Math. Soc.} 31 (1956), 200--205.

%\bibitem{vosper2} G. Vosper, Addendum to "The critical pairs of subsets
%of a group of prime order", {\it J. London Math. Soc.} 31 (1956), 280--282.






\end{thebibliography}
\end{document}